\documentclass[12pt,dvips]{article}
\usepackage{amsmath,amsfonts,amssymb,amsthm,graphicx,verbatim,subfig,bm}
\usepackage[all]{xy}

\ifx\pdftexversion\undefined

\usepackage[a4paper,colorlinks,link
=black,filecolor=black,citecolor=black,urlcolor=black,pdfstartview=FitH]{hyperref}
\else

\usepackage[a4paper,colorlinks,linkcolor=black,filecolor=black,citecolor=black,urlcolor=black,pdfstartview=FitH]{hyperref}
\fi

\usepackage[all]{xy}
\usepackage[a4paper,left=25mm]{geometry}
\ifx\pdftexversion\undefined

\usepackage[a4paper,colorlinks,linkcolor=black,filecolor=black,citecolor=black,urlcolor=black,pdfstartview=FitH]{hyperref}
\fi

\font\sixbb=msbm6
\font\eightbb=msbm8
\font\twelvebb=msbm10 scaled 1095
\newfam\bbfam
\textfont\bbfam=\twelvebb \scriptfont\bbfam=\eightbb
                           \scriptscriptfont\bbfam=\sixbb

\newcommand{\FF}{\mathbb{F}}
\newcommand{\KK}{\mathbb{K}}

\newcommand{\II}{\mathbb{I}}

\newtheorem{theorem}{\bf Theorem}[section]
\newtheorem{claim}[theorem]{\bf Claim}

\newtheorem{proposition}[theorem]{\bf Proposition}
\newtheorem{corollary}[theorem]{\bf Corollary}
\newcommand{\enp}{\begin{flushright} $\Box$ \end{flushright}}
\newcommand{\beq}[0]{\begin{equation}}
\newcommand{\enq}[0]{\end{equation}}

\newcommand{\cb}{{\cal B}}
\newcommand{\cd}{\mathcal{D}}
\newcommand{\cu}{{\cal U}}

\newcommand{\thh}{\tilde{H}}
\newcommand{\supp}{{\rm supp}}
\newcommand{\cf}{{\cal F}}
\newcommand{\cg}{{\cal G}}

\newcommand{\tilb}{\tilde{B}}
\newcommand{\sd}{{\rm sd}}
\newcommand{\sgn}{{\rm sgn}}
\newcommand{\frakg}{\mathfrak{g}}

\newcommand{\frakh}{\mathfrak{h}}

\newcommand{\qbin}[2]{\genfrac{[}{]}{0pt}{}{#1}{#2}}
\newcommand{\spn}{{\rm span} \,}

\newcommand{\mb}[1]{{\mathbf #1}}

\newcommand{\caly}{\mathcal{Y}}
\newcommand{\ccs}{\mathbb{S}}
\newcommand{\epsb}{\pmb{\epsilon}}
\newcommand{\tcv}{\widetilde{c}_{\mb{v}}}

\title{On Lusztig-Dupont Homology \\ of Flag Complexes}
\begin{document}
\author{Roy Meshulam\thanks{Department of Mathematics,
Technion, Haifa 32000, Israel. e-mail:
meshulam@technion.ac.il~. Supported by ISF and GIF grants.} \and Shira Zerbib\thanks{Department of Mathematics,
University of Michigan, Ann Arbor. e-mail: zerbib@umich.edu~. }}
\maketitle
\pagestyle{plain}
\begin{abstract}
Let $V$ be an $n$-dimensional vector space over the finite field $\FF_q$. The spherical building $X_V$ associated with $GL(V)$ is the order complex of the nontrivial linear subspaces of $V$.
Let $\frakg$ be the local coefficient system on $X_V$, whose value on the simplex $\sigma=[V_0 \subset \cdots \subset V_p] \in X_V$ is given by $\frakg(\sigma)=V_0$. The
homology module $\cd^1(V)=\thh_{n-2}(X_V;\frakg)$ plays a key role in Lusztig's seminal work on the discrete series representations of $GL(V)$. Here, some further properties of $\frakg$ and its exterior powers are established.
These include a construction of an explicit basis of $\cd^1(V)$, a computation of the dimension of $\cd^k(V)=\thh_{n-k-1}(X_V;\wedge^k \frakg)$, and the following twisted analogue of a result of Smith and Yoshiara:
For any $1 \leq k \leq n-1$, the minimal support size of a non-zero $(n-k-1)$-cycle in the twisted homology $\thh_{n-k-1}(X_V;\wedge^k \frakg)$
is $\frac{(n-k+2)!}{2}$.
\end{abstract}

\section{Introduction}
\label{s:intro}

Let $q$ be a prime power and let $V$ be an $n$-dimensional vector space over the finite field $\FF_q$. The spherical building associated with $G=GL(V)$ is the order complex $X_V$ of the nontrivial linear subspaces of $V$:
The vertices of $X_V$ are the linear subspaces $0 \neq U \subsetneq V$, and the
$k$-simplices are families of subspaces of the form $\{U_0,\ldots,U_k\}$, where $U_0 \subsetneq \cdots \subsetneq U_k$. The homotopy type of $X_V$ was determined by Solomon and Tits \cite{Solomon} (see also Theorem 4.73 in \cite{AB08}).
\begin{theorem}[Solomon-Tits]
\label{t:build}
$X_V$ is homotopy equivalent to a wedge of $q^{\binom{n}{2}}$ $(n-2)$-spheres. In particular, the reduced homology of $X_V$ with coefficients in a field $\KK$
is given by
$$\dim \thh_i(X_V;\KK)=
\left\{
\begin{array}{ll}
0 & i \neq n-2, \\
q^{\binom{n}{2}} & i=n-2.
\end{array}
\right.~~ $$
\end{theorem}
\ \\ \\
The natural action of $G$ on $X_V$ induces a representation of $G$ on
$\thh_{n-2}(X_V;\KK)$. Viewed as a $G$-module, $\thh_{n-2}(X_V;\KK)$ is
the {\it Steinberg module} of $G$ over $\KK$ (see e.g. section 6.4 in \cite{Smith}).
We recall some facts concerning $X_V$ and the Steinberg module.
For a subset $S\subset V$, let $\langle S \rangle=\spn(S)$ denote the linear span of $S$.
Let $[n]=\{1,\ldots,n\}$.
Let $B=\{v_1,\ldots,v_n\}$ be a basis of $V$ and let $\tilb$ be the set of vertices of $X_V$ given by
$$\tilb=\big\{ \langle v_i:i \in I \rangle~:~ \emptyset \neq I \subsetneq [n]\big\}.$$
The induced subcomplex $X_V[\tilb]$ is the {\it apartment} determined by $B$.
Clearly, $X_V[\tilb]$ is isomorphic to the barycentric subdivision of the boundary of a $(n-1)$-simplex, and thus
$$\thh_{n-2}(X_V[\tilb];\KK) \cong \KK.$$ We next exhibit a generator $z_B$ of $\thh_{n-2}(X_V[\tilb];\KK)$.
For a permutation $\pi$ in the symmetric group $\ccs_n$ and for $1 \leq i \leq n$, let $V_{\pi}(i)=\langle  v_{\pi(1)},\ldots, v_{\pi(i)} \rangle$ and let $\sigma_{\pi}$ be the ordered
$(n-2)$-simplex
$$\sigma_\pi =[V_{\pi}(1) \subset \cdots \subset V_{\pi}(n-1)].$$
Then $z_B=\sum_{\pi \in \ccs_n} {\rm sgn}(\pi) \sigma_{\pi}$ is a generator of $\thh_{n-2}(X_V[\tilb];\KK)$.
The following explicit construction of a basis of $\thh_{n-2}(X_V[\tilb];\KK)$ is due to Solomon \cite{Solomon}
(see also Theorem 4.127 in \cite{AB08}).
\begin{theorem}[Solomon]
\label{t:basis}
Let $\sigma$ be a fixed $(n-2)$-simplex of $X_V$. Then
$$\big\{z_B: B \text{~is~a~basis~of~} V \text{such~that~} \sigma \in X_V[\tilb]\big\}$$
is a basis of $\thh_{n-2}(X_V[\tilb];\KK)$.
\end{theorem}
\noindent
The {\it support} of a $(n-2)$-chain $c=\sum_{\sigma} a_{\sigma} \sigma \in C_{n-2}(X_V;\KK)$ is
$$\supp(c)=\{\sigma:a_{\sigma} \neq 0\}.$$
Clearly, $|\supp(z_B)|=n!$ for any basis $B$ of $V$.
Smith and Yoshiara \cite{SY95} proved that the $z_B$'s are in fact the nontrivial $(n-2)$-cycles of minimal support in $X_V$.
\begin{theorem}[Smith-Yoshiara]
\label{t:wcycle}
$$
\min \big\{ |\supp(z)|: 0 \neq z \in \thh_{n-2}(X_V;\KK) \big\} = n! .
$$
\end{theorem}

In this paper we study analogues of Theorems \ref{t:build}, \ref{t:basis} and \ref{t:wcycle} for the homology of $X_V$ with certain local coefficient systems introduced by Lusztig and Dupont.
We first recall some definitions. Let $X$ be a simplicial complex on a vertex set $S$. Let $\prec$ be an arbitrary fixed linear order on $S$.
For $k \geq -1$ let $X(k)$ denote the set of $k$-dimensional simplices of $X$, and let $X^{(k)}$ denote the $k$-dimensional skeleton of $X$. A simplex
$\sigma \in X(k)$ will be written as $\sigma=[s_1,\ldots,s_{k+1}]$ where $s_1 \prec \cdots \prec s_{k+1}$. The $i$-th face of $\sigma$ as above is the $(k-1)$-simplex $\sigma_i=[s_1,\ldots,s_{i-1},s_{i+1},\ldots,s_{k+1}]$. For a $0$-dimensional simplex $\sigma=
[s_1]$, let $\sigma_1=\emptyset$ be the empty simplex.
A {\it local system} $\cf$ on $X$ is an assignment of an abelian group $\cf(\sigma)$ to each simplex $\sigma \in X$, together with homomorphisms $\rho_{\sigma}^{\tau}:\cf(\tau) \rightarrow \cf(\sigma)$ for each
$\sigma \subset \tau$ satisfying the usual compatibility conditions: $\rho_{\sigma}^{\sigma}={\rm identity}$, and $\rho_{\eta}^{\sigma} \rho_{\sigma}^{\tau}=\rho_{\eta}^{\tau}$ if $\eta \subset \sigma \subset \tau$. A $\cf$-twisted $k$-chain of $X$ is
a formal linear combination $c=\sum_{\sigma \in X(k)} c(\sigma) \sigma$, where $c(\sigma) \in \cf(\sigma)$. Let $C_k(X;\cf)$ denote the group of $\cf$-twisted $k$-chains of $X$.
For $k \geq 0$ define the boundary map
$$\partial_k:C_k(X;\cf) \rightarrow C_{k-1}(X;\cf)$$ by
$$\partial_k \left(\sum_{\sigma \in X(k)} c(\sigma) \sigma\right)=\sum_{\sigma \in X(k)} \sum_{i=1}^{k+1}
(-1)^{i+1} \rho_{\sigma_i}^{\sigma}\big(c(\sigma)\big) \sigma_i.
$$
For $k=-1$ let $\partial_{-1}$ denote the zero map $C_{-1}(X;\cf)=\cf(\emptyset) \rightarrow 0$.
The homology of $X$ with coefficients in $\cf$, denoted by $H_*(X,\cf)$, is the homology of the complex $\oplus_{i \geq 0} C_i(X;\cf)$.
The reduced homology $\thh_*(X,\cf)$ is the homology of $\oplus_{i \geq -1} C_i(X;\cf)$.
Let $X, Y$ be two simplicial complexes and let $f:X \rightarrow Y$ be a simplicial map
such that $\dim f(\sigma)=\dim \sigma$ for all $\sigma \in X$. Let $\cg$ be a local system on $Y$.
The inverse image system $\cf=f^{-1}\cg$ given by $\cf(\sigma)=\cg(f(\sigma))$ is a local system on $X$. The induced mapping on homology is denoted by $f_*:\thh_k(X;\cf) \rightarrow \thh_k(Y;\cg)$.
For further discussion of
local coefficient homology, see e.g. chapter 7 in \cite{Benson91} and chapter 10 in \cite{Smith}.

Lusztig, in his seminal work \cite{L74} on discrete series representations of $GL(V)$, defined and studied the local system $\frakg$ on $X_V$ given by
$\frakg(U_1 \subset \cdots \subset U_{\ell})=U_1$ and $\frakg(\emptyset)=V$, where the connecting homomorphisms $\rho_{\sigma}^{\tau}$'s are the natural inclusion maps. Dupont, in his study of homological approaches to scissors congruences \cite{Dupont82}, extended some of Lusztig's results to the higher exterior powers $\wedge^k \frakg$ over flag complexes of Euclidean spaces.
For $i \geq 0$ let $\thh_i(X_V;\wedge^k \frakg)$ denote the $i$-th homology $\FF_q$-module of the chain complex of $X_V$ with $\wedge^k \frakg$ coefficients. Note that
$C_{-1}(X_V;\wedge^k \frakg)=\wedge^k V$. The following result was proved by Lusztig (Theorem 1.12 in \cite{L74}) for $k=1$, and extended by Dupont (Theorem 3.12 in \cite{Dupont82}) to all $k \geq 1$.
\begin{theorem}[Lusztig, Dupont]
\label{t:ld}
Let $1 \leq k \leq n-1$. Then $\thh_i(X_V;\wedge^k \frakg)=0$ for $i \neq n-k-1$.
\end{theorem}
\noindent
Let $\cd^k(V)=\thh_{n-k-1}(X_V;\wedge^k \frakg)$.
Lusztig (Theorem 1.14 in \cite{L74}) proved that
\begin{equation}
\label{e:dimd1}
\dim \cd^1(V)=\prod_{i=1}^{n-1} (q^i-1).
\end{equation}
The proof of (\ref{e:dimd1}) in \cite{L74} is based on the case $k=1$ of Theorem \ref{t:ld}, combined with an Euler characteristic computation. In Section \ref{s:base1} we describe an explicit
basis of $\cd^1(V)$. This construction may be regarded as a twisted counterpart of Theorem \ref{t:basis}. Concerning the dimension of $\cd^k(V)$ for general $k$, we prove the following extension of Theorem \ref{e:dimd1}.
\begin{theorem}
\label{t:dimdkk}
\begin{equation}
\label{e:dimfor1}
\dim \cd^k(V)=\sum_{1 \leq \alpha_1<\cdots <\alpha_{n-k} \leq n-1} \prod_{j=1}^{n-k}(q^{\alpha_j}-1).
\end{equation}
\end{theorem}
\noindent
Our final result is an analogue of the Smith-Yoshiara Theorem \ref{t:wcycle} for the coefficient system $\wedge^k \frakg$.
\begin{theorem}
\label{t:mindkv}
$$
\min \big\{|\supp(w)|:0 \neq w \in \cd^k(V)\big\}=\frac{(n-k+2)!}{2}.
$$
\end{theorem}
\ \\ \\
The paper is organized as follows. In Section \ref{s:base1} we construct an explicit basis for
$\cd^1(V)$. In Section \ref{s:dimcdk} we use an exact sequence due to Dupont to prove Theorem \ref{t:dimdkk}. In Section \ref{s:mincyc} we recall the Nerve lemma for homology with local coefficients, and obtain a vanishing result for a certain local system on the simplex.
These results are used to prove Theorem \ref{t:mindkv}.
We conclude in Section \ref{s:conc} with some remarks and open problems.

\section{A Basis for $\cd^1(V)$}
\label{s:base1}
In this section we construct an explicit basis for $\cd^1(V)=\thh_{n-2}(X_V;\frakg)$.
Let $V=\FF_q^n$ and let $e_1,\ldots,e_n$ be the standard basis of $V$. For
$a=(a_1,\ldots,a_n), b=(b_1,\ldots,b_n) \in V$ let $a\cdot b$ denote the standard bilinear form
$\sum_{i=1}^n a_i b_i$. For a subset $S\subset V$, let $S^{\perp}=\{u \in V: u\cdot s=0, \text{~for~all~} s \in S\}$.
Let $\prec$ be any linear order on $X_V(0)$ such that $U \prec U'$ if $\dim U < \dim U'$.
Then an $(n-2)$-simplex in $X_V$ is of the form $[U_1,\ldots,U_{n-1}]$, where
$0 \neq U_1 \subsetneq \cdots \subsetneq U_{n-1} \subsetneq V$.

\begin{figure}[h]
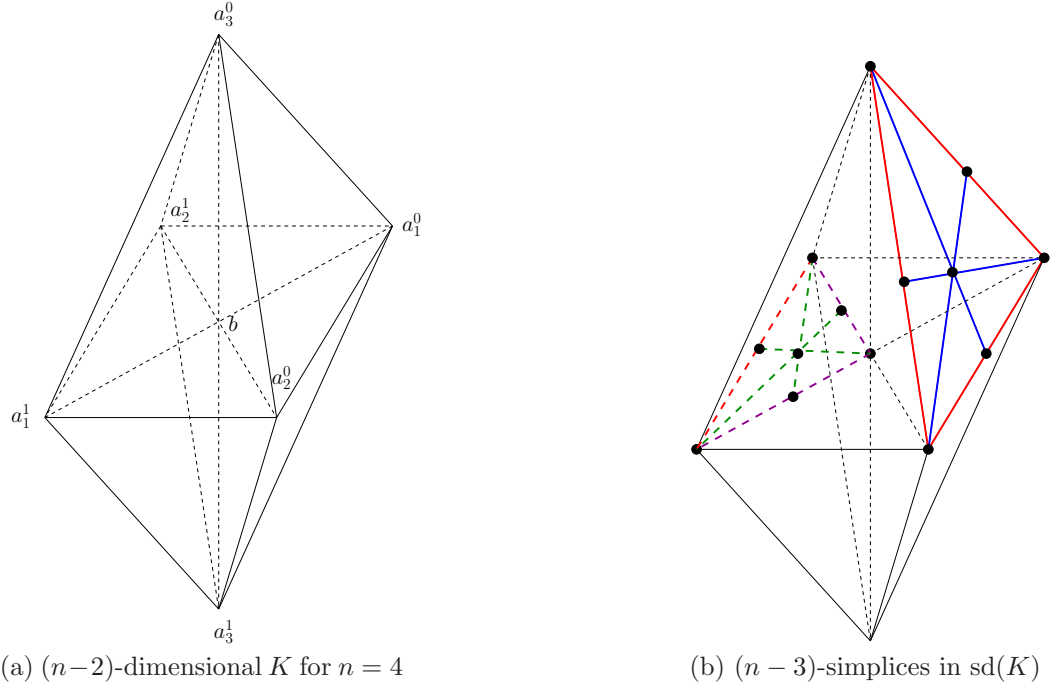

\begin{center}
\subfloat[$(n-2)$-dimensional $K$ for $n=4$]
{\label{fig:kv}
\scalebox{0.4}{\input{octap4.pstex_t}}}
\hspace{100pt}
\subfloat[$(n-3)$-simplices in $\sd(K)$]
{\label{fig:four}
\scalebox{0.4}{\input{octap3.pstex_t}}}
\caption{Four types of $(n-3)$-simplices in $\sd(K)= \supp(c_{\mb{v}})$}
\label{figure3}
\end{center}
\end{figure}

For simplicial complexes $Y$, $Z$ defined on disjoint vertex sets, let
$Y*Z=\{\sigma \cup \tau: \sigma \in Y, \tau \in Z\}$ denote their simplicial join.
Let $a_1^0,a_1^1,\ldots,a_{n-1}^0,a_{n-1}^1,b$ be $2n-1$ distinct elements. Let $M$ denote the octahedral
$(n-2)$-sphere
$$
M=\{a_1^0,a_1^1\}* \cdots * \{a_{n-1}^0,a_{n-1}^1\},
$$
and let $K=M \cup \left(\{b\}*M^{(n-3)}\right)$. See Figure \ref{fig:kv} for a depiction of the
$2$-dimensional complex $K$ when $n=4$.

Choose a linear order $\prec_1$ on the simplices of $K$ such that $\sigma \prec_1 \tau$ if $\dim \sigma > \dim \tau$. The barycentric subdivision of $K$, denoted by $\sd(K)$, is the complex whose vertex set $\sd(K)(0)$ consists of all nonempty simplices of $K$, and whose $k$-simplices (ordered according to $\prec_1$) are
$[\sigma_1,\ldots,\sigma_{k+1}]$ where $\sigma_1 \supsetneq \cdots \supsetneq \sigma_{k+1}$.
For a sequence $\mb{x}=(x_1,\ldots,x_{n-1})$ of distinct vertices of $K$, such that
$\{x_1,\ldots,x_{n-1}\} \in K$, let $S(\mb{x})$ denote the $(n-2)$-simplex of $\sd(K)$
given by
$$
S(\mb{x})=\left[\{x_1,\ldots,x_{n-1}\},\{x_1,\ldots,x_{n-2}\},\ldots,\{x_1\}\right].
$$
For a permutation $\pi$ in the symmetric group $\ccs_{n-1}$ let
$\pi(\mb{x})=(x_{\pi(1)},\ldots,x_{\pi(n-1)})$.
Let $E=\{0,1\}^{n-1}$ and for $1 \leq j \leq n-1$ let
$$E_j=\{(\epsilon_1,\ldots,\epsilon_{n-1}) \in E: \epsilon_j=0\}.$$
For
$\epsb=(\epsilon_1,\ldots,\epsilon_{n-1}) \in E$ and $1 \leq j \leq n-1$
let
$\mb{a}^{\epsb}=(a_1^{\epsilon_1},\ldots,a_{n-1}^{\epsilon_{n-1}})$
and let
$$
\mb{a}^{\epsb,j}=(a_1^{\epsilon_1},\ldots,a_{j-1}^{\epsilon_{j-1}},b,a_{j+1}^{\epsilon_{j+1}},
\ldots,a_{n-1}^{\epsilon_{n-1}}).
$$
Let $T_{q,n}$ denote the set of all sequences $\mb{v}=(v_1,\ldots,v_{n-1}) \in V^{n-1}$ such that
$v_i \in e_i + \langle e_{i+1},\ldots,e_n \rangle$ and $v_i \neq e_i$ for all $1 \leq i \leq n-1$.
Clearly $|T_{q,n}|=\prod_{i=1}^{n-1}(q^i-1)$.
\\
Fix $\mb{v}=(v_1,\ldots,v_{n-1}) \in T_{q,n}$.
For $\epsb=(\epsilon_1,\ldots,\epsilon_{n-1}) \in E$, let $\mb{v}^{\mb{\epsilon}}=(u_1,\ldots,u_{n-1})$, where
$$
u_i=
\left\{
\begin{array}{ll}
e_i & \epsilon_i=0, \\
v_i & \epsilon_i=1.
\end{array}
\right.~~
$$
For $1 \leq j \leq n-1$ let
$\mb{v}^{\epsb,j}=(u_1,\ldots,u_{n-1})$, where
$$u_i=
\left\{
\begin{array}{ll}
e_n & i=j, \\
e_i & i \neq j ~\&~ \epsilon_i=0, \\
v_i & i \neq j ~\&~ \epsilon_i=1.
\end{array}
\right.~~
$$
Define $\theta_{\mb{v}}: K(0) \rightarrow V$
by
$$
\theta_{\mb{v}}(x)= \left\{
\begin{array}{ll}
e_i & x=a_i^0, \\
v_i & x=a_i^1, \\
e_n & x=b,
\end{array}
\right.
$$
and let $f_{\mb{v}}: \sd(K)(0) \rightarrow X_V(0)$
be the map given by
$$
f_{\mb{v}}(\sigma)=\langle \theta_{\mb{v}}(x): x \in \sigma \rangle^{\perp}.
$$
Clearly, $f_{\mb{v}}$ extends to a simplicial map from $\sd(K)$ to $X_V$.
The inverse of $\frakg$ under $f_{\mb{v}}$ is the local system of $\sd(K)$ given by  $\frakh_{\mb{v}}=f_{\mb{v}}^{-1}\frakg$. We next define an element
$$
c_{\mb{v}}=\sum_{F \in \sd(K)(n-2)} c_{\mb{v}}(F) F \in C_{n-2}(\sd(K);\frakh_{\mb{v}}).
$$
For a sequence $\mb{u}=(u_1,\ldots,u_{n-1}) \in V^{n-1}$ of linearly independent vectors in $V$
such that $e_n \not\in \langle u_1,\ldots,u_{n-1} \rangle$, let $w(\mb{u})$ be the unique element $w\in \langle u_1,\ldots,u_{n-1} \rangle^{\perp}$ such that $w \cdot e_n=1$.
For $\epsb=(\epsilon_1,\ldots,\epsilon_{n-1}) \in \{0,1\}^{n-1}$ and $\pi \in \ccs_{n-1}$
let $\chi(\mb{\epsilon},\pi)=(-1)^{\sum_{j=1}^{n-1}\epsilon_j} \sgn(\pi)$.
On an $(n-2)$-simplex $F \in \sd(K)(n-2)$ define
\begin{equation}
\label{e:defc}
c_{\mb{v}}(F)=
\left\{
\begin{array}{ll}
\chi(\epsb,\pi) w(\mb{v}^{\epsb}) &
\epsb \in E,  F=S(\pi(\mb{a}^{\epsb})),\\
\chi(\epsb,\pi) \big(w(\mb{v}^{\epsb+e_j})-w(\mb{v}^{\epsb})\big)
&  \epsb \in E_j, F=S(\pi(\mb{a}^{\epsb,j})), \\
0 & \text{otherwise}.
\end{array}
\right.~~
\end{equation}
\noindent
Note that $c_{\mb{v}}(F) \in \frakh_{\mb{v}}(F)$ for all $F \in \sd(K)(n-2)$.
Indeed, if $F=S(\pi(\mb{a}^{\epsb}))$ then
\begin{equation*}
\begin{split}
c_{\mb{v}}(F)&=\chi(\epsb,\pi)w(\mb{v}^{\epsb}) \in
\langle v_1^{\epsilon_1},\ldots,v_{n-1}^{\epsilon_{n-1}} \rangle^{\perp} \\
&=\frakg(f_{\mb{v}}(F))=\frakh_{\mb{v}}(F).
\end{split}
\end{equation*}
If $F=S(\pi(\mb{a}^{\epsb,j}))$ for $1 \leq j \leq n-1$ and $\epsb \in  E_j$ then
\begin{equation*}
\begin{split}
c_{\mb{v}}(F)&=\chi(\epsb,\pi) \big(w(\mb{v}^{\epsb+e_j})-w(\mb{v}^{\epsb})\big)
\in
\langle v_1^{\epsilon_1},\ldots,v_{j-1}^{\epsilon_{j-1}},e_n,
v_{j+1}^{\epsilon_{j+1}},\ldots,v_{n-1}^{\epsilon_{n-1}} \rangle^{\perp} \\
&=\frakg(f_{\mb{v}}(F))=\frakh(F).
\end{split}
\end{equation*}

\begin{proposition}
\label{p:cv}
$c_{\mb{v}} \in \thh_{n-2}(\sd(K);\frakh_{\mb{v}})$.
\end{proposition}
\noindent
{\bf Proof.}
Let $G \in \sd(K)(n-3)$. We have to show that
$\partial_{n-2} c_{\mb{v}} (G)=0$.
Let $\Gamma(G)$ denote the set of $(n-2)$-simplices in $\sd(K)$ that contain $G$.
For $2 \leq \ell \leq n-1$ let $\eta_{\ell} \in \ccs_{n-1}$ denote the transposition $(n-\ell,n-\ell+1)$. We consider the following four cases according to the type of $G$.
For $n=4$ we depict the various edge types in Figure \ref{fig:four}. Edges colored with blue, green, red and magenta correspond to types 1,2,3 and 4 below respectively.
\begin{enumerate}
\item
$G=S(\pi(\mb{a}^{\epsb}))_{\ell}$ for some $2 \leq \ell \leq n-1$, $\pi \in \ccs_{n-1}$ and
$\epsb \in E$.
\\
Then
$$
\Gamma(G)=\{S(\pi(\mb{a}^{\epsb})),S((\pi \eta_{\ell})(\mb{a}^{\epsb}))\}
$$
As $G$ is the $\ell$-th face of both these simplices, it follows that
\begin{equation*}
\label{e:case1}
\begin{split}
(-1)^{\ell+1}&\partial_{n-2}c_{\mb{v}}(G) =  c_{\mb{v}}(S(\pi(\mb{a}^{\epsb})))+
c_{\mb{v}}(S((\pi \eta_{\ell})(\mb{a}^{\epsb}))) \\
&= \chi(\epsb,\pi) w(\mb{v}^{\epsb})+
\chi(\epsb,\pi \eta_{\ell}) w(\mb{v}^{\epsb})\\
&=\chi(\epsb,\pi) w(\mb{v}^{\epsb})-
\chi(\epsb,\pi) w(\mb{v}^{\epsb})=0.
\end{split}
\end{equation*}

\item
$G=S(\pi(\mb{a}^{\epsb,j}))_{\ell}$ for some $2 \leq \ell \leq n-1$,
$\pi \in \ccs_{n-1}$, $1 \leq j \leq n-1$ and $\epsb \in E_j$.
\\
Then $$\Gamma(G)=\{S(\pi(\mb{a}^{\epsb,j})),
S((\pi \eta_{\ell})(\mb{a}^{\epsb,j}))\}.$$
As $G$ is the $\ell$-th face of both these simplices, it follows that
\begin{equation*}
\label{e:case2}
\begin{split}
(-1)^{\ell+1}&\partial_{n-2}c_{\mb{v}}(G) =  c_{\mb{v}}(S(\pi(\mb{a}^{\epsb,j})))+
c_{\mb{v}}(S((\pi \eta_{\ell})(\mb{a}^{\epsb,j}))) \\
&=\chi(\epsb,\pi) \big(w(\mb{v}^{\epsb+e_j})-w(\mb{v}^{\epsb})\big)+
\chi(\epsb,\pi \eta_{\ell}) \big(w(\mb{v}^{\epsb+e_j})-w(\mb{v}^{\epsb})\big) \\
&=\chi(\epsb,\pi) \big(w(\mb{v}^{\epsb+e_j})-w(\mb{v}^{\epsb})\big)-
\chi(\epsb,\pi) \big(w(\mb{v}^{\epsb+e_j})-w(\mb{v}^{\epsb})\big)
=0.
\end{split}
\end{equation*}

\item
$G=S(\pi(\mb{a}^{\epsb}))_1$ for some $\pi \in \ccs_{n-1}$ and $\epsb \in E_j$,
where $j=\pi(n-1)$.
\\
Then $$\Gamma(G)=\{S(\pi(\mb{a}^{\epsb})), S(\pi(\mb{a}^{\epsb+e_j})),S(\pi(\mb{a}^{\epsb,j}))\}.$$
As $G$ is the $1$-face of each of these simplices, it follows that
\begin{equation*}
\label{e:case3}
\begin{split}
&\partial_{n-2}c_{\mb{v}}(G) =
c_{\mb{v}}\left(S(\pi(\mb{a}^{\epsb}))\right)+c_{\mb{v}}\left( S(\pi(\mb{a}^{\epsb+e_j}))\right)+c_{\mb{v}}\left(S(\pi(\mb{a}^{\epsb,j}))\right) \\
&=\chi(\epsb,\pi) w(\mb{v}^{\epsb})+\chi(\epsb+e_j,\pi) w(\mb{v}^{\epsb+e_j})+\chi(\epsb,\pi) \big(w(\mb{v}^{\epsb+e_j})-w(\mb{v}^{\epsb})\big) \\
&=\chi(\epsb,\pi) \left(w(\mb{v}^{\epsb})- w(\mb{v}^{\epsb+e_j})\right)
+\chi(\epsb,\pi) \big(w(\mb{v}^{\epsb+e_j})-w(\mb{v}^{\epsb})\big)=0.
\end{split}
\end{equation*}

\item
$G=S(\pi(\mb{a}^{\epsb,j}))_1$ for some $\pi \in \ccs_{n-1}$ and $\epsb \in E_j$,
where $j \neq \pi(n-1)$.
\\
Let $j'=\pi(n-1)$ and let $\tau$ denote the transposition $(j,j')$. Since $S(\pi(\mb{a}^{\epsb,j}))_1$ is independent of $\epsilon_{\pi(n-1)}$, we may assume that
$\epsilon_{j'}=\epsilon_{\pi(n-1)}=0$. Then:
$$\Gamma(G)=\{S(\pi(\mb{a}^{\epsb,j})),S(\pi(\mb{a}^{\epsb+e_{j'},j})), S((\tau \pi)(\mb{a}^{\epsb,j'})),S((\tau \pi)(\mb{a}^{\epsb+e_j,j'}))\}.$$

As $G$ is the $1$-face of each of these simplices, it follows that
\begin{equation*}
\label{e:case4}
\begin{split}
\partial_{n-2}c_{\mb{v}}(G) &=
 c_{\mb{v}}(S(\pi(\mb{a}^{\epsb,j})))
+c_{\mb{v}}(S(\pi(\mb{a}^{\epsb+e_{j'},j}))) \\
&~~~~~ +c_{\mb{v}}(S((\tau \pi)(\mb{a}^{\epsb,j'})))
+c_{\mb{v}}(S((\tau \pi)(\mb{a}^{\epsb+e_j,j'}))) \\
&=\chi(\epsb,\pi)(w(\mb{v}^{\epsb+e_j})-w(\mb{v}^{\epsb}))+
\chi(\epsb+e_{j'},\pi)
(w(\mb{v}^{\epsb+e_{j'}+e_j})-w(\mb{v}^{\epsb+e_{j'}})) \\
&~~~~~+\chi(\epsb,\tau\pi)(w(\mb{v}^{\epsb+e_{j'}})-w(\mb{v}^{\epsb}))+
\chi(\epsb+e_{j},\tau\pi)
(w(\mb{v}^{\epsb+e_{j}+e_{j'}})-w(\mb{v}^{\epsb+e_{j}})) \\
&=\chi(\epsb,\pi)[(w(\mb{v}^{\epsb+e_j})-w(\mb{v}^{\epsb}))
-(w(\mb{v}^{\epsb+e_{j'}+e_j})-w(\mb{v}^{\epsb+e_{j'}})) \\
&~~~~~-(w(\mb{v}^{\epsb+e_{j'}})-w(\mb{v}^{\epsb}))+
(w(\mb{v}^{\epsb+e_{j'}+e_{j}})-w(\mb{v}^{\epsb+e_{j}}))]=0.
\end{split}
\end{equation*}
\end{enumerate}
We have thus shown that $c_{\mb{v}} \in \thh_{n-2}(\sd(K);\frakh)$.
{\enp}
\noindent
Proposition \ref{p:cv} implies that $\tcv=(f_{\mb{v}})_*c_{\mb{v}} \in \thh_{n-2}(X_V;\frakg)$.
\begin{theorem}
\label{t:dcyc}
The family $\{\tcv: \mb{v} \in T_{q,n}\}$ is a basis of $\cd^1(V)=\thh_{n-2}(X_V;\frakg)$.
\end{theorem}
\noindent
{\bf Proof.} Let $\mb{v} \in T_{q,n}$.
Let $R(\mb{v}) \in X_V(n-2)$ be the $(n-2)$-simplex
$$
R(\mb{v})=
[\langle v_1,\ldots,v_{n-1}\rangle^{\perp}, \langle v_1,\ldots,v_{n-2}\rangle^{\perp},\ldots,
\langle v_1,v_2\rangle^{\perp},
\langle v_1 \rangle^{\perp}].
$$
Let $\mb{1}=(1,\ldots,1) \in E$. It is straightforward to check that
$F=S(\mb{a}^{\mb{1}})$ is the unique $(n-2)$-simplex in $\sd(K)$ such that
$f_{\mb{v}}(F)=R(\mb{v})$. It follows that
$$\tcv (R(\mb{v}))=c_{\mb{v}}(S(\mb{a}^{\mb{1}}))=(-1)^{n-1} w(\mb{v}).$$
On the other hand, if $\mb{v} \neq \mb{v'} \in T_{q,n}$, then
$R(\mb{v'}) \not\in f_{\mb{v}}(\sd(K))$ and so $\tcv (R(\mb{v'}))=0$.
It follows that the $(n-2)$-cycles $\{\tcv: \mb{v} \in T_{q,n}\}$
are linearly independent in $\cd^1(V)$. As $|T_{q,n}|=\prod_{i=1}^{n-1}(q^i-1)=\dim\cd^1(V)$, this completes the proof of Theorem \ref{t:dcyc}.
{\enp}
\noindent
{\bf Example:} Let $n=3$ and let
$$\mb{v}=(v_1,v_2)=\left( (1,r,s),(0,1,t) \right) \in T_{q,3}.$$
Figure \ref{fig:snz} depicts the cycle $c_{\mb{v}} \in H_1(\sd(K);\frakh)$.
Black vertices correspond to vertices of $K$ and red vertices correspond to edges of $K$.
The values of $c_{\mb{v}}$ are indicated on the edges of the diagram. For example,
let $\epsb=(1,1)$ and $\pi=(1,2)$. Then
$F=S(\pi(\mb{a}^{\epsb}))=[\{a_2^1,a_1^1\},\{a_2^1\}]$, and
$$c_{\mb{v}}(F) = \chi(\epsb,\pi) w\big((v_1,v_2)\big)=-w\big((v_1,v_2)\big)=(s-rt,t,-1).$$
Similarly, if $j=1$, $\epsb=(0,1) \in E_1$ and $\pi=(1,2)$, then
$F=S(\pi(\mb{a}^{\epsb,j}))=[\{a_2^1,b\},\{a_2^1\}]$ and
\begin{equation*}
\begin{split}
c_{\mb{v}}(F) &= \chi(\epsb,\pi) \left( w\big((v_1,v_2)\big)-w\big((e_1,v_2)\big) \right) \\
&=(rt-s,-t,1)-(0,-t,1)=(rt-s,0,0).
\end{split}
\end{equation*}

\begin{figure}[h!]
\begin{center}
  \scalebox{0.4}{\input{sdk2.pstex_t}}
  \caption{The cycle $c_{\mb{v}}$ for $\mb{v}=(v_1,v_2)=\left( (1,r,s),(0,1,t)\right)$.}
  \label{fig:snz}
\end{center}
\end{figure}
\noindent
Figures \ref{fig:genb} and \ref{fig:degb} depict
the $1$-cycle $\tcv \in H_1(X_V;\frakg)$.
Here, the black vertices correspond to $2$-dimensional subspaces of $V$. The red vertices and their labels correspond to $1$-dimensional subspaces and their generating vectors. Figure \ref{fig:genb} depicts the generic case when $rst(rt-s) \neq 0$. The labels of the left most $6$ red points together with the $\pm$ signs, indicate the values of $\tcv$ on the incident edges. The remaining three values of $\tcv$ are indicated on the edges incident with the vertex corresponding to the line spanned by $(1,0,0)$. Figure \ref{fig:degb} similarly depicts the case $s=0$. Note that in both cases, the simplicial map $f_{\mb{v}}:\sd(K) \rightarrow X_V$ is not injective.

\begin{figure}[h!]
\begin{center}
  \scalebox{0.4}{\input{genbase.pstex_t}}
  \caption{The cycle $\tcv$ for a generic $\mb{v}=\left( (1,r,s),(0,1,t)\right)$.}
  \label{fig:genb}
\end{center}
\end{figure}

\begin{figure}[h!]
\begin{center}
  \scalebox{0.4}{\input{degbaser.pstex_t}}
  \caption{The cycle $\tcv$ for $\mb{v}=(v_1,v_2)=\left( (1,r,0),(0,1,t)\right)$.}
  \label{fig:degb}
\end{center}
\end{figure}

\section{The Dimension of $\cd^k(V)$}
\label{s:dimcdk}
\noindent
{\bf Proof of Theorem \ref{t:dimdkk}:}
For an $\FF_q$-space $W$ let $\text{St}(W)=\thh_{\dim W -2}(X_W;\FF_q)$ denote the Steinberg module of $W$ over $\FF_q$.
Recall that  $\dim \text{St}(W)=q^{\binom{\dim W}{2}}$ by Theorem \ref{t:build}.
Let $G_j(V)$ denote the family of all $j$-dimensional linear subspaces of $V$.
The following result is due to Dupont (Proposition 5.38 in \cite{Dupont82}).
\begin{theorem}[Dupont]
\label{dups1}
There is an exact sequence
\begin{equation*}
\label{dup1}
\begin{split}
0 \rightarrow \cd^k(V) &\rightarrow \bigoplus_{U_k \in G_k(V)} \wedge^{\, k} U_k \otimes {\rm St}(V/U_k) \rightarrow \bigoplus_{U_{k+1} \in G_{k+1}(V)} \wedge^k U_{k+1} \otimes {\rm St}(V/U_k) \rightarrow
\\
\ldots &\rightarrow \bigoplus_{U_{n-2} \in G_{n-2}(V)} \wedge^k U_{n-2}\otimes {\rm St}(V/U_{n-2})
\rightarrow \bigoplus_{U_{n-1} \in G_{n-1}(V)} \wedge^k U_{n-1} \rightarrow \wedge^k V \rightarrow 0.
\end{split}
\end{equation*}
\end{theorem}
\noindent
Writing $\qbin{n}{j}_q$ for the $q$-binomial coefficient,
Theorem \ref{dups1} implies that
\begin{equation}
\label{dimdk}
\dim \cd^k(V)=\sum_{j=k}^n (-1)^{j-k} \binom{j}{k} q^{\binom{n-j}{2}} \qbin{n}{j}_q.
\end{equation}
\noindent
By the $q$-binomial theorem (see e.g. (1.87) in \cite{Stanley})
\begin{equation}
\label{qbinid}
\prod_{j=0}^{n-1} (1+q^j \lambda)=\sum_{j=0}^n q^{\binom{j}{2}} \qbin{n}{j}_q \lambda^j.
\end{equation}
Substituting $\lambda=-t^{-1}$ in (\ref{qbinid}) and multiplying by $t^n$ it follows that
\begin{equation}
\label{st1}
\prod_{j=0}^{n-1} (t-q^j)=\sum_{j=0}^n (-1)^j q^{\binom{j}{2}} \qbin{n}{j}_q t^{n-j}.
\end{equation}
Differentiating (\ref{st1}) $k$ times and multiplying by $\frac{(-1)^{n-k}}{k!}$ we obtain
\begin{equation}
\label{difk}
\begin{split}
& \prod_{j=0}^{n-1} (q^j-t)\sum_{0 \leq \alpha_0<\cdots < \alpha_{k-1} \leq n-1} \prod_{\ell=0}^{k-1} \frac{1}{q^{\alpha_{\ell}}-t} \\
&= \sum_{j=0}^n (-1)^{n-k+j} \binom{n-j}{k} q^{\binom{j}{2}} \qbin{n}{j}_q t^{n-j-k}  \\
&= \sum_{j=0}^n (-1)^{j-k} \binom{j}{k} q^{\binom{n-j}{2}} \qbin{n}{j}_q t^{j-k}.
\end{split}
\end{equation}
Substituting $t=1$ in (\ref{difk}) and using (\ref{dimdk}) we obtain (\ref{e:dimfor1}).
{\enp}

\subsection{A Basis for $\cd^{n-1}(V)$}
\label{subs:cdno}
In this subsection we describe an explicit basis for $\cd^{n-1}(V)=\thh_{0}(X_V;\wedge^{n-1} \frakg)$. We first recall some facts concerning the exterior algebra $\wedge V$.
Let $V=\FF_q^n$. Using the notation of Section \ref{s:base1}, recall that $e_1,\ldots,e_n$ are the unit vectors in $V$, and $a \cdot b$ denotes the standard symmetric bilinear form on $V$.
Let $\mb{e}=e_1 \wedge \cdots \wedge e_n \in \wedge^n V$.
The induced bilinear form on $\wedge^p V$ is given by
$$(u_1 \wedge \cdots \wedge u_p) \cdot (v_1\wedge \cdots \wedge v_p)= \det
\big(u_i \cdot v_j\big)_{i,j=1}^p.$$
The {\it star operator} $*:\wedge^{n-k} V \rightarrow \wedge^{k}V$ is the unique linear map
that satisfies $$(*\alpha) \cdot \beta= \mb{e} \cdot (\alpha \wedge \beta)$$ for any
$\alpha \in \wedge^{n-k} V, \beta \in \wedge^{k}V$.
\begin{claim}
\label{c:dualb}
Let $v_1,\ldots,v_{n-k}$ be linearly independent vectors in $V$ and let
$M= \langle v_1,\ldots,v_{n-k} \rangle^{\perp}$.
Then $$0 \neq *(v_1 \wedge \cdots \wedge v_{n-k})\in \wedge^k M.$$
\end{claim}
\noindent
{\bf Proof.} Extend $\{v_i\}_{i=1}^{n-k}$ to a basis $\{v_i\}_{i=1}^n$ of $V$, and let
$\{w_j\}_{j=1}^n$ be the dual basis, i.e. $v_i \cdot w_j=\delta_{i,j}$. Then
$M=\langle w_{n-k+1},\ldots,w_n \rangle$.
For a subset $L=\{i_1,\ldots,i_{\ell}\} \in \binom{[n]}{\ell}$ such that
$1 \leq i_1 < \cdots < i_{\ell} \leq n$ let $v_L=v_{i_1} \wedge \cdots \wedge v_{i_{\ell}}$
and $w_L=w_{i_1} \wedge \cdots \wedge w_{i_{\ell}}$. If $L,L' \in \binom{[n]}{\ell}$ then
$v_L\cdot w_{L'}=\delta_{L,L'}$.
\\
Let $I_0=\{1,\ldots,n-k\}$, $J_0=\{n-k+1,\ldots,n\}$,
and let $*v_{I_0}=\sum_{|J|=k} \lambda_J w_J$. Then for any $J' \in \binom{[n]}{k}$
\begin{equation}
\label{e:dotp1}
*v_{I_0}\cdot v_{J'}=\sum_{|J|=k} \lambda_J w_J \cdot v_{J'}=\lambda_{J'}.
\end{equation}
On the other hand
\begin{equation}
\label{e:dotp2}
\begin{split}
*v_{I_0} \cdot v_{J'}&=\mb{e} \cdot (v_{I_0} \wedge v_{J'}) \\
&=
\left\{
\begin{array}{ll}
\det (v_1,\ldots,v_n) & J'=J_0, \\
0 & J' \neq J_0.
\end{array}
\right.~~
\end{split}
\end{equation}
Combining (\ref{e:dotp1}) and (\ref{e:dotp2}), it follows that
$0 \neq *v_{I_0}= \det (v_1,\ldots,v_n) w_{J_0} \in \wedge^k M$.
{\enp}
\noindent
We proceed to construct a basis of $\cd^{n-1}(V)=\thh_{0}(X_V;\wedge^{n-1} \frakg)$.
Note that if $u \in V$, then by Claim \ref{c:dualb}, $(*u)u^{\perp} \in C_0(X_V;\wedge^{n-1}\frakg)$.  For any $1 \leq i \leq n$ let
$$z_{u,i}=(*e_i)e_i^{\perp}+(*u)u^{\perp}-(*(u+e_i))(u+e_i)^{\perp} \in
C_0(X_V;\wedge^{n-1}\frakg).$$
Then $$\partial_0(z_{u,i})=*e_i+*u-*(u+e_i)=*(e_i+u-(u+e_i))=0$$
and therefore $z_{u,i} \in \cd^{n-1}(V)$. For $2 \leq i \leq n$ let
$R_i= \left(\FF_q^{i-1}\setminus \{0\}\right) \times \{0\}^{n-i+1}$.
\begin{claim}
\label{c:bdnmo}
\begin{equation}
\label{e:bdnmo}
\cb=\big\{z_{u,i}: 2 \leq i \leq n~,~u \in R_i\big\}
\end{equation}
is a basis of $\cd^{n-1}(V)$.
\end{claim}
\noindent
{\bf Proof.}
By Theorem \ref{t:dimdkk}
$$\dim \cd^{n-1}(V)=\sum_{i=2}^n (q^{i-1}-1)=\sum_{i=2}^n |R_i|=|\cb|.$$
It therefore suffices to show that the elements of $\cb$ are linearly independent.
This in turn follows from the fact that for any $2  \leq j \leq n$  and
$v \in R_j$, it holds that $(v+e_j)^{\perp} \in \supp(z_{v,j})$,
but $(v+e_j)^{\perp} \not\in \supp(z_{u,i})$
for any $(u,i) \neq (v,j)$ such that $2 \leq i \leq j$ and $u \in R_i$.
{\enp}

\section{Minimal Cycles in $\cd^k(V)$}
\label{s:mincyc}
\noindent
In this section we prove Theorem \ref{t:mindkv}.
The upper bound follows from a construction of certain explicit $(n-k-1)$-cycles of $\cd^k(V)$ given in Subsection \ref{sub:upperb}. The lower bound is established in Subsection \ref{sub:lowerb}.

\subsection{The Upper Bound}
\label{sub:upperb}

Let $1 \leq k \leq n-1$ and let $m=n-k+2$. Let $\mb{u}=(u_1,\ldots,u_m) \in V^m$ be an ordered $m$-tuple of vectors in $V$ whose only linear dependence is $\sum_{i=1}^m u_i=0$.
Let $\II_{m-2,m}$ denote the family of injective functions
$\pi:[n-k]=[m-2] \rightarrow [m]$. For $\pi \in \II_{m-2,m}$ let $T(\mb{u},\pi)$ be the
$(n-k-1)$-simplex given by
$$
T(\mb{u},\pi)=[\langle u_{\pi(1)},\ldots,u_{\pi(n-k)}\rangle^{\perp} \subset
\cdots  \subset  \langle u_{\pi(1)}\rangle^{\perp}].
$$
Let $\gamma_{\mb{u}} \in C_{n-k-1}(X_V;\wedge^k\frakg)$ be the chain whose value on an $(n-k-1)$-simplex $F$ is given by
\begin{equation}
\label{e:defgam}
\gamma_{\mb{u}}(F)=
\left\{
\begin{array}{ll}
*\left(u_{\pi(1)} \wedge \cdots \wedge u_{\pi(n-k)}\right) & F=T(\mb{u},\pi), \\
0 & \text{otherwise}.
\end{array}
\right.~~
\end{equation}

\begin{proposition}
\label{p:cycwk}
$\gamma_{\mb{u}} \in \cd^k(V)$.
\end{proposition}
\noindent
{\bf Proof.} Let $G$ be an $(n-k-2)$-simplex in $X_V$. Let $\Gamma_{\mb{u}}(G)$ denote the set of $(n-k-1)$-simplices in $\supp(\gamma_{\mb{u}})$ that contain $G$.
For $2 \leq \ell \leq n-k$ let $\eta_{\ell} \in \ccs_{n-k-2}$ denote the transposition $(n-k-\ell+1,n-k-\ell+2)$. We consider the following two cases:
\begin{enumerate}
\item
$G=T(\mb{u},\pi)_{\ell}$ for some $2 \leq \ell \leq n-k$ and $\pi \in \II_{m-2,m}$.
\\
Then
$$
\Gamma_{\mb{u}}(G)=\left\{ T(\mb{u},\pi),T(\mb{u},\pi\eta_{\ell})\right\}.
$$
As $G$ is the $\ell$-th face of both these simplices, it follows that
\begin{equation*}
\label{e:case1}
\begin{split}
(-1)^{\ell+1}&\partial_{n-k-1} \gamma_{\mb{u}}(G) =
\gamma_{\mb{u}}(T(\mb{u},\pi))+\gamma_{\mb{u}}(T(\mb{u},\pi\eta_{\ell})) \\
&=*\left(u_{\pi(1)} \wedge \cdots \wedge u_{\pi(n-k-\ell+1)} \wedge u_{\pi(n-k-\ell+2)} \wedge \cdots \wedge u_{\pi(n-k)}\right) \\
&~~+*\left(u_{\pi(1)} \wedge \cdots \wedge u_{\pi(n-k-\ell+2)} \wedge u_{\pi(n-k-\ell+1)} \wedge \cdots \wedge u_{\pi(n-k)}\right)=0.
\end{split}
\end{equation*}
\item
$G=T(\mb{u},\pi)_{1}$ for some $\pi \in \II_{m-2,m}$.
\\
Let $[m] \setminus \pi([m-3])= \{\alpha_1,\alpha_2,\alpha_3\}$.
For $i=1,2,3$ define $\pi_i \in \II_{m-2,m}$ by
$$
\pi_i(j)=\left\{
\begin{array}{ll}
\pi(j)  & 1 \leq j \leq n-k-1, \\
\alpha_i & j=n-k.
\end{array}
\right.~~
$$
Then
$$
\Gamma_{\mb{u}}(G)=\left\{ T(\mb{u},\pi_1),T(\mb{u},\pi_2),T(\mb{u},\pi_3)\right\}.
$$
As $G$ is the $1$-th face of these three simplices, it follows that
\begin{equation*}
\label{e:case1}
\begin{split}
&\partial_{n-k-1} \gamma_{\mb{u}}(G) =
\sum_{i=1}^3 \gamma_{\mb{u}}(T(\mb{u},\pi_i)) \\
&= \sum_{i=1}^3 *\big(u_{\pi(1)} \wedge \cdots \wedge u_{\pi(n-k-1)} \wedge u_{\alpha_i}\big) \\
&= *\big(u_{\pi(1)} \wedge \cdots \wedge u_{\pi(n-k-1)} \wedge (\sum_{i=1}^3 u_{\alpha_i})\big) \\
&=*\big(u_{\pi(1)} \wedge \cdots \wedge u_{\pi(n-k-1)} \wedge (\sum_{j=1}^m u_{j})\big)=0.
\end{split}
\end{equation*}
\end{enumerate}
We have thus shown that $\gamma_{\mb{u}} \in \cd^k(V)$.
{\enp}
\begin{corollary}
\label{c:lbnd}
\begin{equation*}
\label{e:lbnd}
\begin{split}
\min &\big\{|\supp(w)|:0 \neq w \in \cd^k(V)\big\} \leq |\supp(\gamma_{\mb{u}})| \\
&=|\II_{m-2,m}|=\frac{(n-k+2)!}{2}.
\end{split}
\end{equation*}
\end{corollary}
\noindent
{\bf Example:} Let $n=3, k=1$. A minimal twisted $1$-cycle in $\cd^1(X_V)$ is depicted in Figure \ref{fig:degb}.

\subsection{The Lower Bound}
\label{sub:lowerb}

In preparation for the proof of the lower bound in Theorem \ref{t:mindkv}, we first recall a twisted version of the nerve lemma. Let $\cf$ be a local system on a finite simplicial complex $Y$, and
let $\caly=\{Y_i\}_{i=1}^m$ be a family of subcomplexes of $Y$ such that $Y=\bigcup_{i=1}^m Y_i$.
The {\it nerve} of the cover $\caly$ is the simplicial complex $N=N(\caly)$ on the vertex $[m]=\{1,\ldots,m\}$, whose simplices are the subsets $\tau \subset [m]$ such that $Y_{\tau}:=\bigcap_{i \in \tau} Y_i \neq \emptyset$. For $j \geq 1$ let $N_j(\cf)$ be the local system on $N$ given by
$N_j(\cf)(\tau)=H_j(Y_{\tau};\cf)$. The following result is twisted version of the Mayer-Vietoris spectral sequence (see e.g. \cite{L74}).
\begin{proposition}
\label{mv}
There exists a spectral sequence $\{E_{p,q}^r\}$ converging to $H_{*}(Y;\cf)$
such that $E_{p,q}^1=\bigoplus_{\sigma \in N(p)} H_q(Y_{\sigma};\cf)$ and $E_{p,q}^2=H_p(N;N_q(\cf))$.
\end{proposition}
\noindent
The Nerve Lemma is the following
\begin{corollary}
\label{nerve1}
Suppose that $H_q(Y_{\sigma};\cf)=0$ for all $q \geq 1$ and $\sigma \in N(p)$ such that $p+q \leq t$.
Then $H_p(Y;\cf)\cong H_p(N;N_0(\cf))$ for all $0 \leq p \leq t$.
\end{corollary}
\noindent
We will also need a simple observation concerning a certain twisted homology of the simplex.
Let $r \geq 2$ and let $W_1,\ldots,W_r$ be arbitrary linear subspaces of a finite
dimensional vector space $W$ over a field $\KK$.
Let $\Delta_{r-1}$ denote the simplex on the vertex set $[r]$, and
let $\cg$ be the local system on $\Delta_{r-1}$ given by
$$
\cg(\sigma)=
\left\{
\begin{array}{ll}
\bigcap_{i \in \sigma} W_i & \emptyset \neq \sigma \in \Delta_{r-1}, \\
W & \sigma=\emptyset,
\end{array}
\right.~~
$$
with the natural inclusion maps.
\begin{proposition}
\label{p:gfact}
$\thh_{k}(\Delta_{r-1};\cg)=0$ for $k \geq r-2$.
\end{proposition}
\noindent
{\bf Proof:} Using the natural order on $\{1,\ldots,r\}$, the top dimensional simplex in $\Delta_{r-1}$ is
$\tau=[1,2,\cdots,r]$, and its
$i$-th face is $\tau_i=[1,\ldots,i-1,i+1,\ldots,r]$. For $1 \leq i < j \leq r$ let
$$\tau_{i,j}=[1,\ldots,i-1,i+1,\ldots,j-1,j+1,\ldots,r].$$ Then
$$C_{r-1}(\Delta_{r-1},\cg)=\left\{w \tau: w \in \bigcap_{i=1}^r W_i\right\}$$ and
$$C_{r-2}(\Delta_{r-1},\cg)= \left\{\sum_{i=1}^r w_i \tau_i: w_i \in \bigcap_{\ell \in \tau_i} W_{\ell}\right\}.$$
The boundary map $\partial_{r-1}:C_{r-1}(\Delta_{r-1};\cg) \rightarrow C_{r-2}(\Delta_{r-1};\cg)$ is  given by
\begin{equation}
\label{e:br1}
\partial_{r-1}(w\tau)=\sum_{i=1}^r (-1)^{i+1} w\tau_i.
\end{equation}
Note that for $1 \leq i \leq r$ and $1 \leq j \leq r-1$, the $j$-th face of $\tau_i$ is
\begin{equation*}
(\tau_i)_j= \left\{
\begin{array}{ll}
\tau_{j,i} & 1 \leq j <i \leq r, \\
\tau_{i,j+1} & 1 \leq i \leq j \leq r-1.
\end{array}
\right.
\end{equation*}
It follows that the boundary map $\partial_{r-2}:C_{r-2}(\Delta_{r-1};\cg) \rightarrow C_{r-3}(\Delta_{r-1};\cg)$ is given by
\begin{equation}
\label{e:br2}
\begin{split}
\partial_{r-2} \big(\sum_{i=1}^r w_i \tau_i\big)&=\sum_{i=1}^r \sum_{j=1}^{r-1} (-1)^{j+1} w_i (\tau_i)_j \\
&=\sum_{i=1}^r \sum_{j=1}^{i-1} (-1)^{j+1} w_i \tau_{j,i} +
\sum_{i=1}^r \sum_{j=i}^{r-1} (-1)^{j+1} w_i \tau_{i,j+1} \\
&=\sum_{j=1}^r \sum_{i=1}^{j-1} (-1)^{i+1} w_j \tau_{i,j} +
\sum_{i=1}^r \sum_{j=i+1}^{r} (-1)^{j} w_i \tau_{i,j} \\
&=\sum_{1 \leq i <j \leq r} \left( (-1)^{i+1} w_j+(-1)^j w_i\right) \tau_{i,j}.
\end{split}
\end{equation}
Eq. (\ref{e:br1}) implies that $\thh_{r-1}(\Delta_{r-1};\cg)=0$. Next
let $c=\sum_{i=1}^r w_i \tau_i \in \ker \partial_{r-2}$ be a $\cg$-twisted $(r-2)$-cycle of $\Delta_{r-1}$.
It follows by (\ref{e:br2})
that $w_j=(-1)^{j+1}w_1$ for all $1 \leq j \leq r$. Therefore $w_1 \in \bigcap_{i=1}^r W_i$ and hence $w_1 \tau \in C_{r-1}(X;\cg)$. Eq. (\ref{e:br1}) then implies that $\partial_{r-1}(w_1 \tau)=c$. Thus $\thh_{r-2}(\Delta_{r-1};\cg)=0$.
{\enp}

\noindent
{\bf Proof of the lower bound in Theorem \ref{t:mindkv}}.
We argue by induction on $n-k$. For the induction basis $k=n-1$, we have to show that
if $0 \neq z \in \cd^{n-1}(V)=\thh_0(X_V;\wedge^{n-1} \frakg)$, then
$|\supp(z)| \geq \frac{(n-k+2)!}{2}=3$.
Suppose for contradiction that $|\supp(z)| < 3$.
Then $z=(*u) u^{\perp}+(*v) v^{\perp}$ for some $u,v \in V$.
As $$0=\partial_0 z=(*u)+(*v)=*(u+v),$$ it follows that $u+v=0$ and hence $z=0$, a contradiction.
For the induction step, assume that $n-k \geq 2$ and let
$$0 \neq z =\sum_{z \in X_V(n-k-1)} z(\tau) \tau \in H_{n-k-1}(X_V;\wedge^k \frakg)=Z_{n-k-1}(X_V;\wedge^k \frakg).$$
Let $\supp(z)=\{\tau_1,\ldots,\tau_s\} \in X_V(n-k-1)$ and write
$$\tau_i=[V_k(i),\ldots,V_{n-1}(i)],$$ where $\dim V_j(i)=j$ for all $1 \leq i \leq s$ and $k \leq j \leq n-1$.
Let
$$\{V_{n-1}(i):1 \leq i \leq s\}=\{U_1,\ldots,U_{r}\},$$
where the $U_i$'s are distinct $(n-1)$-dimensional subspaces. Let
$\cu_i=\{U:0 \neq U \subset U_i\}$
and let $Y_i=X_V[\cu_i]$. Let $Y=\cup_{i=1}^{r} Y_i$ then clearly $z \in Z_{n-k-1}(Y;\wedge^k \frakg)$.
Let $N$ be the nerve of the cover $\{Y_i\}_{i=1}^{r}$ of $Y$.
For $\sigma \subset [r]$ let $U_{\sigma}=\cap_{i \in \sigma} U_i$ and $Y_{\sigma}=\cap_{i \in \sigma} Y_i$.
If $\sigma \in N$ then $Y_{\sigma}$ is the order complex
of the poset $P_{\sigma}=\{W:0 \neq W \subset U_{\sigma}\}$. As $P_{\sigma}$ has a unique maximal element $U_{\sigma}$
it follows (see e.g. Lemma 1.4 in \cite{L74}) that
\begin{equation}
\label{nqq}
N_q(\wedge^k \frakg)(\sigma)= H_q(Y_{\sigma};\wedge^k \frakg)=
\left\{
\begin{array}{ll}
\wedge^k U_{\sigma} & q=0, \\
0 & q>0.
\end{array}
\right.~~
\end{equation}
Write
$$\cf(\sigma)=N_0(\wedge^k \frakg)(\sigma)=\wedge^k U_{\sigma}.$$
Eq. (\ref{nqq}) and
Corollary \ref{nerve1} imply that for all $p \geq 0$
\begin{equation}
\label{e:ncons}
H_p(Y;\wedge^k \frakg) \cong H_p(N;\cf).
\end{equation}
\begin{proposition}
\label{boundl}
$r \geq n-k+2$.
\end{proposition}
\noindent
{\bf Proof:}
Suppose to the contrary that $r \leq n-k+1$. Then $\Delta_{r-1}^{(r-2)} \subset N \subset \Delta_{r-1}$.
For $1 \leq i \leq r$ let
$W_i= \wedge^k U_i \subset \wedge^k V$. Let $\cg$ be the local system on $\Delta_{r-1}$ given
by $\cg(\sigma)=\cap_{i \in \sigma} W_i=\wedge^k U_{\sigma}$. Then $\cg(\sigma)=\cf(\sigma)$ if
$\sigma \in N$ and $\cg(\sigma)=0$ otherwise.
Hence $H_*(\Delta_{r-1};\cg)=H_*(N;\cf)$.
As $n-k-1 \geq r-2$, it follows by
combining (\ref{e:ncons}) and Proposition \ref{p:gfact} that
$$H_{n-k-1}(Y;\wedge^k \frakg) \cong H_{n-k-1}(N;\cf)=H_{n-k-1}(\Delta_{r-1};\cg)=0,$$
in contradiction with the assumption that $z$ is a nonzero element of $H_{n-k-1}(Y;\wedge^k \frakg)$.
{\enp}
\noindent
We now conclude the proof of Theorem \ref{t:mindkv}.
For $1 \leq j \leq r$ define
$z_j \in C_{n-k-2}(X_{U_j};\wedge^k \frakg)$ as follows. For an $(n-k-2)$-simplex
$F=[V_k,\ldots,V_{n-2}] \in X_{U_j}(n-k-2)$ let
$z_j(F) =z([V_k,\ldots,V_{n-2},U_j])$.
Then  $\partial_{n-k-2}z_j=0$. Indeed, suppose that
$$[V_k,\ldots,V_{i-1},V_{i+1},\ldots,V_{n-2}] \in X_{U_j}(n-k-3),$$ where
$\dim V_{\ell}=\ell$ for $i \neq \ell \in \{k,\ldots,n-2\}$. Then:
\begin{equation*}
\begin{split}
&\partial_{n-k-2} z_j\left([V_k,\ldots,V_{i-1},V_{i+1},\ldots,V_{n-2}]\right) \\
&= (-1)^{i+k} \sum_{V_{i-1} \subset V_i \subset V_{i+1}}
z_j\left([V_k,\ldots,V_{i-1},V_i,V_{i+1},\ldots,V_{n-2}]\right) \\
&=(-1)^{i+k} \sum_{V_{i-1} \subset V_i \subset V_{i+1}}
z\left([V_k,\ldots,V_{i-1},V_i,V_{i+1},\ldots,V_{n-2},U_{j}]\right) \\
&= \partial_{n-k-1} z\left([V_k,\ldots,V_{i-1},V_{i+1},\ldots,V_{n-2},U_j]\right)=0.
\end{split}
\end{equation*}
As  $0 \neq z_j \in H_{n-k-2}(X_{U_j};\wedge^k \frakg)$, it follows
by induction that $|\supp(z_j)| \geq \frac{(n-k+1)!}{2}$.  Therefore by Proposition \ref{boundl}
$$|\supp(z)|=\sum_{j=1}^{r} |\supp(z_j)| \geq (n-k+2)\frac{(n-k+1)!}{2}=\frac{(n-k+2)!}{2}.$$
{\enp}

\section{Concluding Remarks}
\label{s:conc}
In this paper we studied some aspects of the twisted homology modules
$\cd^k(V)=\thh_{n-k-1}(X_V;\wedge^k \frakg)$. Our results suggest several problems and directions for further research:
\begin{itemize}
\item
In Sections \ref{s:base1} and \ref{subs:cdno} we described explicit bases for $\cd^1(V)=\thh_{n-2}(X_V;\frakg)$ and for $\cd^{n-1}(V)=\thh_{0}(X_V;\wedge^{n-1} \frakg)$.
It would be interesting to obtain analogous constructions for other $\cd^k(V)$'s.

\item
The Nerve Lemma argument used in the proof of Theorem \ref{t:mindkv}
can be adapted to give a simple alternative proof of the Smith-Yoshiara Theorem \ref{t:wcycle}.
We hope that this approach can also be useful for the study of minimal cycles of local
systems over other highly symmetric complexes.

\item
The Smith-Yoshiara Theorem \ref{t:wcycle} and its counterpart for the local system
$\wedge^k \frakg$, Theorem \ref{t:mindkv}, show that the linear codes that arise from (twisted) homology of $X_V$ have small distance relative to their length, and are therefore far from good codes. On the other hand, it is known (see \cite{ALLM13}) that for fixed integers $n \geq 2$ and $K>0$  there is a constant $\lambda=\lambda(n,K)>0$,
such that for sufficiently large $N$ there exists a complex $X_N \subset \Delta_{N-1}^{(n)}$
whose number of $n$-faces satisfies $f_n(X_N)= K \binom{N}{n}$, and such that $|\supp(z)| \geq \lambda \binom{N}{n}$ for all
$0 \neq z \in C=H_n(X_N;\FF_2)$.
In particular, the rate $r(C)$ and relative distance $\delta(C)$ of $C$ satisfy
$$r(C)=\frac{\dim C}{f_n(X_N)} \geq \frac{K-1}{K}$$
and
$$\delta(C)= \frac{\min\{|\supp(c)|:0 \neq c \in C\}}{f_n(X_N)}\geq \frac{\lambda}{K}.$$
It would be interesting to give explicit constructions of simplicial complexes that give rise to homological codes with similar parameters.

\end{itemize}

\end{document}